\input amstex
\documentstyle{amsppt}
\magnification=\magstep1 \NoRunningHeads
\topmatter
\title Weak mixing for nonsingular Bernoulli\\  actions of countable amenable groups
\endtitle

\author
Alexandre I. Danilenko 
\endauthor

\email
alexandre.danilenko@gmail.com
\endemail

\address
 Institute for Low Temperature Physics
\& Engineering of National Academy of Sciences of Ukraine, 47 Nauky Ave.,
 Kharkiv, 61103, UKRAINE
\endaddress
\email alexandre.danilenko\@gmail.com
\endemail

\abstract
Let $G$ be an amenable   discrete countable infinite group,
 $A$  a finite set, and $(\mu_g)_{g\in G}$ a family of probability measures on $A$ such that $\inf_{g\in G}\min_{a\in A}\mu_g(a)>0$.
It is shown   (among other results)  that if  the Bernoulli shiftwise action of $G$ on the infinite product  space $\bigotimes_{g\in G}(A,\mu_g)$
is nonsingular and conservative then it is weakly mixing.
This answers in positive a question  by Z.~Kosloff who proved recently that the conservative Bernoulli $\Bbb Z^d$-actions are ergodic.
As a byproduct, we  prove a weak version of the pointwise ratio ergodic theorem for nonsingular actions of $G$.
\endabstract

\endtopmatter

\document

\head 0. Introduction
\endhead

Let $G$ be an amenable  discrete  infinite  countable  group.
We call a nonsingular $G$-action $Q=(Q_g)_{g\in G}$ on a $\sigma$-finite measure space    {\it weakly mixing} (see \cite{GlWe}, \cite{DaSi} and references therein) if for each ergodic measure preserving $G$-action $R=(R_g)_{g\in G}$
on a standard probability space, the product $G$-action $(Q_g\times R_g)_{g\in G}$ is ergodic.
Of course, every weakly mixing action is ergodic.

Given a countable set $A$ and a family $(\mu_g)_{g\in G}$ of probability measures on $A$,
we set $X:= A^G$ and $\mu:=\bigotimes_{g\in G}\mu_g$.
Throughout  this section we denote by $T=(T_g)_{g\in G}$ the left shiftwise action of $G$ on $X$, i.e.
$(T_gx)_h:=x_{g^{-1}h}$ for all $x=(x_h)_{h\in G}\in X$ and $g,h\in G$.
The dynamical system $(X,\mu, T)$ is called  a {\it nonsingular Bernoulli $G$-action} if 
$\mu\circ T_g\sim\mu$ for all $g\in G$.
If $\mu_g=\mu_h$ for all $g,h\in G$ then $T$ preserves $\mu$ and the dynamical properties of 
probability preserving Bernoulli actions are well understood (see \cite{OrWe}).% for amenable $G$ and more recent  \cite{Bo}, \cite{Se} with references therein for general $G$).

On the other hand, the purely nonsingular case is considerably less studied even for $G=\Bbb Z$.
Some of nonsingular Bernoulli shifts can be nonconservative (see examples in \cite{Ha} and \cite{DaLe}).
Krengel  constructed in \cite{Kr} the first  conservative nonsingular Bernoulli shift %$T$ which is not of Krieger's type $II_1$, i.e. $T$
which  does not admit an equivalent invariant probability measure (see also \cite{Ha} and \cite{Ko1} for further refinements of his result).
It was assumed in all those papers that $\mu$ is {\it semistationary}, i.e. there is $n\in\Bbb Z$ such that either $\mu_n=\mu_{n-1}=\mu_{n-2}=~\cdots$ or $\mu_n=\mu_{n+1}=\mu_{n+2}=~\cdots$.
In the recent works \cite{Ko3} and \cite{DaLe}, the semistationary nonsingular Bernoulli shifts were  studied  in depth  in  the framework of theory of nonsingular endomorphisms.
In particular, it was shown that every such shift is either dissipative or weakly mixing, it  possesses    the nonsingular property $K$ (in the sense of \cite{SiTh}) and    Krieger's type of it (in the conservative case) is either $II_1$ or  $III_1$.
A number of explicit examples of type $III_1$ Bernoulli transformations with various weak mixing properties were constructed by Vaes and Wahl in \cite{VaWa}.
They also showed that each infinite countable  amenable group (in fact, a group from a much larger class of countable groups with nontrivial
first $L^2$-cohomology)
 has   Bernoulli actions of
type $III_1$.

 %Moreover, if the shift is of type $III_1$ then it and the Maharam extension of it posesses the nonsingular property $K$.
\comment
The  semistationarity implies that
 \roster
 \item"$(\bullet)$"
 the Maharam extension of the natural extension of a 1-sided shift is the natural extension of the Maharam extension of the shift.
 \endroster
This fact  was used essentially in \cite{Ko3} and \cite{DaLe}.
\endcomment

Unfortunately, the aforementioned approach to semistationary shifts ($\Bbb Z$-actions) does not work with the  shifts $T$ which do not have an equivalent semistationary  measure. 
Such shifts  remain almost unstudied so far.
It is known that they are of zero type \cite{Ko2}.
Some progress was achieved recently by Kosloff \cite{Ko4} who proved that
\roster
\item"$(\circ)$" if $A$ is finite,  $\inf_{n\in\Bbb Z}\min_{a\in A}\mu_n(a)>0$ and $T$ is conservative then $T$ is ergodic.
\endroster 
We refine and extend this as follows:

\proclaim{Theorem 0.1 \rom{(see Corollary 3.2, Theorem 3.5 below)}}
\roster
\item"$(i)$"
  For a finite $A$, if $\inf_{n\in\Bbb N}\min_{a\in A}\mu_n(a)>0$ and $T$ is conservative then $T$ is weakly mixing.
\item"$(ii)$"
If $\# A=2$, 
$\inf_{n\in\Bbb N}\min_{a\in A}\Big|\log\Big(\frac{\mu_n(a)}{\mu_{n+1}(a)}\Big)\Big|<\infty$ and $T$ is conservative then $T$ is weakly mixing.
\item"$(iii)$" Under condition $(i)$ or $(ii)$, if $T\times\cdots\times T(p\text{\ times})$ is conservative for some $p>1$ then  $T\times\cdots\times T(p\text{\ times})$ is weakly mixing.
\endroster
\endproclaim

We note that the condition on $(\mu_n)_{n=1}^\infty$ in  $(ii)$ is weaker than the condition on $(\mu_n)_{n=1}^\infty$ in  $(i)$.
Hence $(i)$ and $(\circ)$ follow from  $(ii)$ in case where $\# A=2$.
We also construct an explicit example of a weakly mixing nonsingular Bernoulli shift whose quasiinvariant measure is not equivalent to any semistationary one (see Example~3.3 below).

The proof of $(\circ)$ in \cite{Ko4} is based heavily on application  of the Hurewicz nonsingular ergodic theorem and the maximal inequality.
That proof  is valid also  for the nonsingular Bernoulli actions of $\Bbb Z^d$, $d<\infty$, and, apparently, some Heisenberg groups, however it  does not extend to actions of  groups for which
the Hurewicz theorem fails,  say $\bigoplus_{n=1}^\infty\Bbb Z$ or $\Bbb Q$  (see a discussion in \cite{Ho}). 
In view of that Kosloff asks \cite{Ko4, Problem from \S5.3}: 
\roster
\item""
does $(\circ)$ extend to   nonsingular Bernoulli actions of arbitrary countable amenable group actions?
\endroster
We answer in positive by demonstrating a stronger theorem.

\proclaim{Theorem 0.2 \rom{(see Corollary 2.4 below)}}
 Let $A$ be finite and let  $G$ be an amenable  discrete countable infinite group.
  \roster
 \item"$(i)$"
If $\inf_{g\in G}\min_{a\in A}\mu_g(a)>0$ and $T=(T_g)_{g\in G}$ is conservative then $T$ is weakly mixing.
 \item"$(ii)$"
 Under the condition  of $(i)$,
 if the ``diagonal'' $G$-action
 $T\times\cdots\times T(p\text{\ times})$ is conservative for some $p>1$ then  it is weakly mixing.
\endroster
\endproclaim

\comment

In this connection, we note that if $G$ has vanishing first $L^2$-cohomology, i.e. $H^1(G,\ell^2(G))=\{0\}$, then every conservative nonsingular Bernoulli $G$-action on $A^G$ admits an equivalent $G$-invariant Bernoulli measure $\lambda^{\otimes G}$, where $\lambda$ is a probability on $A$ \cite{VaWa, Theorem~3.1}. 
Since each probability preserving Bernoulli action is mixing,
 the assertion of Theorem~0.2  for such a $G$ follows immediately from the fact that
 mixing implies  weak mixing \cite{Sc}.

\endcomment

In fact, we  deduce Theorem 0.2 from the following more general result (cf. \cite{Ko4, Theorem~2}).

\proclaim{Theorem 0.3 \rom{(see Theorem 2.3 below)}} Let $(X,d)$ be a Polish ultrametric space,  $\Cal R$ a
countable  Borel equivalence relation on $X$, 
 $\mu$  a probability Borel measure on $X$ and
 $Q=(Q_g)_{g\in G}$ a conservative nonsingular $G$-action on $(X,\mu)$.
 If there is a $\mu$-conull subset $X_0\subset X$ such that
 \roster
 \item"---"
 $(Q_g\times Q_g)(\Cal R\cap (X_0\times X_0))=\Cal R\cap(X_0\times X_0)$  for all $g\in G$, 
  \item"---"
 if $(Q_gx,x)\in\Cal R$ for some $x\in X_0$ and $g\in G$ then $g=1$,
  \item"---"
the restriction of $Q$ to the $\sigma$-algebra  $\text{\rom{Inv}}(\Cal R)$ of $\Cal R$-saturated subsets is ergodic,
 \item"---"
$\lim_{g\to\infty}d(Q_gx,Q_gy)=0$ for all $(x,y)\in\Cal R\cap(X_0\times X_0)$,
 \item"---"
there is a Borel map $\alpha:\Cal R\to [1,+\infty)$  with
$$
\alpha(x,y)^{-1}<\frac{d\mu\circ Q_g}{d\mu}(x)/\frac{d\mu\circ Q_g}{d\mu}(y)<\alpha(x,y)
$$
for all $(x,y)\in\Cal R\cap(X_0\times X_0)$ and each $g\in G$.
\endroster
then $Q$ is ergodic.
If, moreover, $\text{\rom{Inv}}(\Cal R)$ is trivial $\pmod \mu$ then $Q$ is weakly mixing.
\endproclaim

In contrast to \cite{Ko4}, we do not use any form of the maximal inequality  in the proof of Theorem~0.3. 
Our main tools are the measurable orbit theory and the following version
of the pointwise ratio ergodic theorem for nonsingular group actions which is of independent interest.

\proclaim{Theorem 0.4 \rom{(see Appendix)}}
Let $G$ be amenable.
Then for each  nonsingular action $T=(T_g)_{g\in G}$ on a standard probability space $(X,\mu)$ and a countable subset $\Cal L$ of $L^1(X,\mu)$, there is 
 a F{\o}lner sequence  $(F_n)_{n=1}^\infty$ in $G$ 
such that $\bigcup_{n=1}^\infty F_n=G$
and for each $f\in\Cal L$,
 there exists
$$
\lim_{n\to\infty}\frac{\sum_{g\in F_n}f( T_gx)\frac{d\mu\circ T_g}{d\mu}(x)}{\sum_{g\in F_n}\frac{d\mu\circ T_g}{d\mu}(x)}=E(f\mid\goth T)(x)
$$
at a.e. $x$, where $\goth T$ denotes the $\sigma$-algebra of $T$-invariant subsets in $X$.
\endproclaim

Nonsingular Markov shifts were also considered in \cite{Ko4}.
It was shown there that under a natural boundedness condition, if the shifts are conservative then they are ergodic. 
We deduce from Theorem~0.3 that they are weakly mixing.

The outline of the paper is as follows.
In \S1 we briefly remind some basic concepts from the theory of measured equivalence relations.
In  \S2 we prove Theorems~0.2 and 0.3.
In the final \S3 we consider the case $G=\Bbb Z$ in more detail.
We first extend slightly Theorems~0.2  in this particular case (see Theorem~3.1) and then prove
Theorem~0.1.
The example of  a conservative nonsingular Bernoulli shift whose quasiinvariant measure is not equivalent to semistationary  one is also constructed there.
A remark on nonsingular Markov shifts completes \S3.
Appendix contains a proof of Theorem~0.4.

{\it Acknowledgements.} I am grateful to M.~Lema{\'n}czyk  and S. Vaes for useful discussions and for finding  some  gaps in  earlier versions of the paper.
I also thank Z.~Kosloff for useful remarks.

\head 1. Countable equivalence relations on measure spaces
\endhead

Let $(X,\goth B)$ be a standard Borel space.
A Borel equivalence relation $\Cal R\subset X\times X$ is called {\it countable} if for each $x\in X$, the $\Cal R$-class $\Cal R(x):=\{y\in X\mid (x,y)\in\Cal R\}$ of $x$ is countable.
If $\Cal R(x)$ is finite for each $x$ then $\Cal R$ is called {\it finite}.

Given a subset $B\in\goth B$, we denote by $\Cal R(B)$ the $\Cal R$-saturation of $B$,
i.e. $\Cal R(B):=\bigcup_{x\in B}\Cal R(x)$. 
Of course,  $\Cal R(B)\in\goth B$.
A subset $B\in\goth B$ is called {\it $\Cal R$-saturated} (or {\it $\Cal R$-invariant}) if $B=\Cal R(B)$.
We denote the collection of all $\Cal R$-saturated subsets by Inv$(\Cal R)$.
It is a sub-$\sigma$-algebra of $\goth B$.
If  $T$ is a Borel bijection $T$ of $X$  and
 $T\Cal R(x)=\Cal R(Tx)$ at each $x\in X$, we say that $T$ {\it normalizes $\Cal R$}.
In this case the restriction of $T$ to Inv$(\Cal R)$  is  well defined. 

Let $A$ and $B$ be two countable sets,  $\# A\ge 2$, $\# D=\infty$, and $X=A^D$.
We define {\it the tail equivalence relation} on $X$ by setting that a point $x=(x_d)_{d\in D}\in X$ is equivalent to a point $y=(y_d)_{d\in D}\in X$ if there is a finite subset $D_0\subset D$ such that $x_d=y_d$ whenever $d\not\in D_0$.

Let $\mu$ be a probability measure on $(X,\goth B)$.
We say that $\Cal R$ is {\it $\mu$-nonsingular} if
$\mu(\Cal R(B))=0$ whenever $\mu(B)=0$.
If %$\Cal R$ is $\mu$-nonsingular and 
each subset from Inv$(\Cal R)$ 
is either $\mu$-null or $\mu$-conull then $\Cal R$ is called {\it $\mu$-ergodic}.
If %$\Cal R$ is $\mu$-nonsingular and  
$\Cal R(x)\cap B\ne\emptyset$ for each subset $B\in\goth B$  for a.e. $x\in B$
then $\Cal R$ is called {\it $\mu$-conservative}.
If $\Cal R$ is $\mu$-nonsingular and
 there is an increasing sequence $\Cal R_1\subset\Cal R_2\subset\cdots$ of finite equivalence relations $\Cal R_n$, $n\in\Bbb N$, such that $\Cal R(x)=\bigcup_{n=1}^\infty\Cal R_n(x)$ for a.e. $x\in X$ then $\Cal R$ is called {\it $\mu$-hyperfinite} and the sequence $(\Cal R_n)_{n=1}^\infty$ is called {\it a filtration of $\Cal R$}.
 A Borel equivalence subrelation of a $\mu$-nonsingular equivalence relation
 is  $\mu$-nonsingular itself.
 \comment
 It is straightforward to verify that for each $\mu$-nonsingular equivalence relation $\Cal R$, there exists a hyperfinite equivalence subrelation $\Cal S\subset\Cal R$ such that
 Inv$(\Cal R)=\text{Inv}(\Cal S)$ (mod 0).
 Hence $\Cal R$ is $\mu$-ergodic if and only if $\Cal S$ is $\mu$-ergodic.
 Moreover,
  $\Cal R$ is 
  $\mu$-conservative  if and only if  $\Cal S$ is $\mu$-conservative.
  Indeed, consider the disintegration of $\mu$ with respect to Inv$(\Cal R)$. 
  Then $\Cal R$ is $\mu$-conservative if and only if  a.e. conditional measure in this disintegration is either non-atomic or supported on a  finite $\Cal R$-class.
Thus the property of $\Cal R$ to be $\mu$-conservative is determined by Inv$(\Cal R)$.

\endcomment
  
For each countable $\Cal R$, there is a countable group $\Gamma$ of Borel isomorphisms of $X$ such that $\Cal R$ is the orbit equivalence relation of $\Gamma$ \cite{FeMo}.
Such a group $\Gamma$ is not unique.
However $\Cal R$ is $\mu$-nonsingular, $\mu$-ergodic or $\mu$-conservative if and only if 
$\Gamma$ is $\mu$-nonsingular, $\mu$-ergodic or $\mu$-conservative respectively.

If $\Cal R$ is $\mu$-nonsingular then there exists a Borel map $\Delta_{\Cal R,\mu}:\Cal R\to\Bbb R_+^*$ such that
$$
\Delta_{\Cal R,\mu}(x,y)=\Delta_{\Cal R,\mu}(y,z)\Delta_{\Cal R,\mu}(z,x)\quad
\text{for all $(x,y),(y,z)\in\Cal R$}
$$
 and $\Delta_{\Cal R,\mu}(x,\gamma x)=\frac{d\mu\circ\gamma}{d\mu}(x)$ at a.e. $x$ for each $\mu$-nonsingular invertible Borel transformation $\gamma$ of $X$ with $(x,\gamma x)\in \Cal R$ for a.e. $x$.
 We call $\Delta_{\Cal R,\mu}$ {\it the Radon-Nikodym cocycle} of the pair $(\Cal R,\mu)$.
 A $\mu$-nonsingular $\Cal R$ is $\mu$-conservative if and only if $\sum_{y\in\Cal R(x)}\Delta_{\Cal R,\mu}(x,y)=+\infty$ at a.e. $x$.
 If $\Cal S$ is a subrelation of $\Cal R$ then $\Delta_{\Cal S,\mu}=\Delta_{\Cal R,\mu}\restriction\Cal S$.
 If $\Cal R$ is finite and $f\in L^1(X,\mu)$ then the mathematical expectation of $f$ with respect to Inv$(\Cal R)$ at a point $x\in X$ is the ratio
 $$
 \frac{\sum_{y\in\Cal R(x)}f(y)\Delta_{\Cal R,\mu}(x,y)}{\sum_{y\in\Cal R(x)}\Delta_{\Cal R,\mu}(x,y)}.
 $$
For more information about measured equivalence relations (and measurable orbit theory) we refer to \cite{DaSi} and references therein.

\head 2. Proof of  Theorems 0.3 and 0.2.
\endhead

Let $X$ be a Polish 0-dimensional space, $d$ an ultrametric on $X$ compatible with the topology, and   $\Cal R$  a  countable Borel equivalence relation on $X$. 
%and let $\Delta_\Cal R:\Cal R\to\Bbb R$ stand for  the Radon-Nikodym cocycle of $\Cal R$.
Let $G$ be an amenable  discrete infinite countable   group and  let $T=(T_g)_{g\in G}$  be a Borel free action of $G$ on $X$ that  {\it normalizes} $\Cal R$, i.e.  $T_g$ normalizes $\Cal R$
for each $g\in G$. 
%Suppose that $T$ is $\mu$-nonsingular and  
We recall that given   a probability measure $\mu$ on $X$, $T$ is called
 {\it strictly $\Cal R$-outer} (mod $\mu$),  if  $T_gx\not\in\Cal R(x)$ for each $g\in G\setminus\{1\}$  at a.e. $x\in X$.
 Of course, a strictly $\Cal R$-outer  $G$-action is free (mod $\mu$).

\definition{Definition 2.1} \roster
\item"(i)"
We say that $T$  {\it squashes} $\Cal R$ if there is a $\mu$-conull subset $X_0\subset X$ such that $\lim_{g\to\infty}d(T_gx,T_gy)=0$ for each $(x,y)\in\Cal R\cap(X_0\times X_0)$.
\item"(ii)"
We say that $T$  is $\Cal R$-{\it bounded} if
 there are a Borel map $\alpha:\Cal R\to [1,+\infty)$ and a $\mu$-conull subset $X_0\subset X$ such that
$$
\alpha(x,y)^{-1}\le\frac{d\mu\circ T_g}{d\mu}(x)/\frac{d\mu\circ T_g}{d\mu}(y)\le\alpha(x,y)
$$
for all $(x,y)\in\Cal R\cap(X_0\times X_0)$ and each $g\in G$.
\endroster
\enddefinition

\remark{Remark \rom{2.2}} It is easy to verify that if $\Cal R$ is $\mu$-nonsingular then
$T$ is $\Cal R$-bounded if and only if
there are a Borel map $\beta:\Cal R\to [1,+\infty)$ and a $\mu$-conull subset $Y_0\subset X$ such that
$$
\beta(x,y)^{-1}\le\Delta_{\Cal R,\mu}(T_gx,T_gy)\le\beta(x,y)
$$
for all $(x,y)\in\Cal R\cap(Y_0\times Y_0)$ and each $g\in G$.
\endremark

We now state and prove the main result of this section.

\proclaim{Theorem 2.3} Let $\mu$ be a probability Borel measure on $X$ such that
 $T=(T_g)_{g\in G}$ is a conservative nonsingular  $G$-action on $(X,\mu)$.
Let $T$ normalize  a countable  Borel equivalence relation
 $\Cal R$ on $X$ and let the restriction of $T$ to \rom{Inv}$(\Cal R)$ be $\mu$-ergodic
 \footnote{The latter is equivalent to the fact that the smallest equivalence relation  containing $\Cal R$ and the graphs of all $T_g$, $g\in G$, is $\mu$-ergodic.}.  
 If $T$ is strictly $\Cal R$-outer, $\Cal R$-bounded and squashes $\Cal R$ then $T$ is ergodic.
If, moreover, $\Cal R$ is $\mu$-ergodic then $T$ is weakly mixing.
\endproclaim

\demo{Proof} 
Suppose that $T$ is not ergodic.
Denote by $\goth F$ the $\sigma$-algebra of Borel $T$-invariant subsets.
Then there is a subset $B\in \goth F$ with $0<\mu(B)<1$.
Let $E(.\mid\goth F)$ stand for the mathematical expectation with respect to $\goth F$.
Select a sequence $(A_n)_{n=1}^\infty$ of  open subsets in $X$ such that 
$$
\frac 1{n^3}\ge \|1_B-1_{A_n}\|_1\ge\|1_B-E(1_{A_n}\mid\goth F)\|_1
$$
for each $n$. 
Then by the Markov inequality,
$$
\mu\bigg(\bigg\{x\in X\,\bigg|\, |1_B(x)-E(1_{A_n}\mid\goth F)(x)|>\frac 1n\bigg\}\bigg)\le \frac1{n^2}.
$$
Hence the Borel-Cantelli lemma yields that for a.e. $x\in X$, 
$$
|1_B(x)-E(1_{A_n}\mid\goth F)(x)|\le \frac 1n\quad\text{ eventually in $n$.}\tag2-1
$$
 Since $d$ is an ultrametric, we may assume without loss of generality that
 for each $n>0$, there is  $\epsilon_n>0$ such that for each $x\in A_n$, the ball centered  at $x$ and of radius $\epsilon_n$ is a subset of $A_n$.
 Denote by $\Cal T$ the $T$-orbit equivalence relation on $(X,\mu)$.
 It is countable, $\mu$-nonsingular, $\mu$-conservative and $\mu$-hyperfinite. 
 \comment

Select a hyperfinite equivalence subrelation   $\Cal S$ of $\Cal T$ such that 
Inv$(\Cal S)=\text{Inv}(\Cal T)$ (mod 0).
It follows that $\Cal S$ is $\mu$-nonsingular and $\mu$-conservative.
Fix a filtration $\Cal S_1\subset \Cal S_2\subset\cdots$ of $\Cal S$
with finite equivalence subrelations.
Then $\text{Inv}(\Cal S_1)\supset \text{Inv}(\Cal S_2) \supset\cdots$ and $\bigwedge_{k=1}^\infty\text{Inv}(\Cal S_k)=\text{Inv}(\Cal S)=\text{Inv}(\Cal T)=\goth F$.
It follows from the martingale convergence theorem
that
$E(1_{A_n}\mid\text{Inv}(\Cal S_k))\to E(1_{A_n}\mid\goth F)$ as $k\to\infty $ almost everywhere for each $n>0$.
For $x\in X$ and $k>0$, we denote by $G_{x,k}$ a finite subset of $G$ such that
$\Cal S_k(x)=\{T_gx\mid g\in G_{x,k}\}$.
Then we have that $G_{x,1}\subset G_{x,2}\subset\cdots$, $\bigcup_{m>0}\{T_gx\mid g\in G_{x,m}\}=\Cal S(x)$
and

\endcomment
By Theorem~A.1 (see Appendix below), there is an increasing sequence
$H_1\subset H_2\subset\cdots$ of finite subsets  in $G$ such that
$\bigcup_{k=1}^\infty H_k=G$ and
$$
\aligned
E(1_{A_n}\mid\goth F)(x)%&=\lim_{k\to+\infty}\frac{\sum_{y\in\Cal S_k(x)}1_{A_n}(y)\Delta_{\Cal S,\mu}(x,y)}{\sum_{y\in\Cal S_k(x)}\Delta_{\Cal S,\mu}(x,y)}\\
&=\lim_{k\to+\infty}\frac{\sum_{g\in H_k}1_{A_n}(T_gx)\frac{d\mu\circ T_g}{d\mu}(x)}{\sum_{g\in H_k}\frac{d\mu\circ T_g}{d\mu}(x)}
\endaligned
\tag2-2
$$
for each $n$ at a.e. $x\in X$.
Since $T$ normalizes $\Cal R$, it follows  that $\Cal R(B)\in\goth F$.
Therefore, by the condition of the theorem, $\mu(\Cal R(B))=1$ and thus
$\mu( \Cal R(B)\setminus B)>0$.
Hence there exist a Borel subset $B_0\subset B$, a $\mu$-nonsingular Borel one-to-one map $\gamma:B_0\to X\setminus B$ and  $C>0$ such that  $\mu(B_0)>0$, $(x,\gamma x)\in\Cal R$ and $\alpha(x,\gamma x)<C$  for  all $x\in B_0$.
Since $T$ squashes $\Cal R$, for a.e. $x\in X$ and each $n>0$, there is  a finite subset $F_n(x)$  in $G$
such that $d(T_gx,T_g\gamma x)<\epsilon_n$ whenever $g\not\in F_n(x)$.
According to our choice  of $\epsilon_n$, if $g\not\in F_n(x)$ then $T_gx\in A_n$ if and only if $T_g\gamma x\in A_n$.
Since $T$ is $\Cal R$-bounded,
 we obtain  that
$$
\frac 1{\alpha(x,\gamma x)} \frac{d\mu\circ T_g}{d\mu}(\gamma x)\le\frac{d\mu\circ T_g}{d\mu}(x)\le{\alpha(x,\gamma x)} \frac{d\mu\circ T_g}{d\mu}(\gamma x)\tag2-3
$$
at a.e. $x\in B_0$.
Since
$T$ is conservative,
$$
\infty=%\sum_{y\in\Cal S(x)}\Delta_{\Cal S,\mu}(x,y)=
\sum_{g\in G}\frac{d\mu\circ T_g}{d\mu}(x)=\lim_{k\to\infty}\sum_{g\in H_{k}\setminus F_n(x)}\frac{d\mu\circ T_g}{d\mu}(x)\tag2-4
$$
at a.e. $x$ and each $n>0$.
Substituting \thetag{2-3} into \thetag{2-2}  and using \thetag{2-4} we obtain that
$$
E(1_{A_n}\mid\goth F)(x)\le\alpha(x,\gamma x)^2E(1_{A_n}\mid\goth F)(\gamma x)\le
C^2E(1_{A_n}\mid\goth F)(\gamma x)
$$
for all $n>0$ and a.e. $x\in B_0$.
Passing to the limit   ($n\to\infty$) and applying \thetag{2-2} we obtain that $1\le 0$, which is a nonsense.
This proves the first claim of the theorem.

To prove the second one, take 
 an ergodic probability preserving $G$-action $R=(R_g)_{g\in G}$.
 Since every Borel  $G$-action admits a countable generating partition  \cite{JaKeLo},
we
can realize $R$ as a left shiftwise action  on the space $\Bbb N^G$ which is 0-dimensional Polish in the usual infinite product topology.
We define an ultrametric $d'$ on $\Bbb N^G$ by setting $d'(a,b):=\sum_{h\in G}\epsilon_h(1-\delta(a_h,b_h))$, where $\delta(.,.)$ is the Kronecker delta, $\epsilon_h>0$ for each $h\in G$ and $\sum_{h\in G}\epsilon_h<\infty$.
Denote by $\Cal T$ the tail equivalence relation on $\Bbb N^G$.
Then $R$ normalizes and squashes  $\Cal T$.
The product space $X\times \Bbb N^G$ is 0-dimentional Polish.
Its topology is compatible with an ultrametric.
%Denote by $\Cal T$ the tail equivalence relation on $\Bbb N^G$.
%Of course, $R$ normalizes $\Cal T$.
%Hence
The product $G$-action $T\times R:=(T_g\times R_g)_{g\in G}$ normalizes and squashes the product  equivalence relation $\Cal R\times\Cal T$ on $X\times \Bbb N^G$.
Moreover, $T\times R$ is strictly $(\Cal R\times\Cal T)$-outer because $T$ is strictly $\Cal R$-outer.
Since
$$
\text{Inv}(\Cal R\times\Cal T)=\text{Inv}(\Cal R)\otimes\text{Inv}(\Cal T)
$$
and $\text{Inv}(\Cal R)$ is trivial mod $\mu$,
 the restriction of $T\times R$ to $\text{Inv}(\Cal S\times\Cal T)$
is isomorphic to the restriction of $R$ to $\text{Inv}(\Cal T)$, which is ergodic because $R$ is ergodic on the entire Borel $\sigma$-algebra on $\Bbb N^G$.
It is easy to verify that 
since $R$ is probability preserving,   $T\times R$ is $(\Cal R\times\Cal T)$-bounded.
Since $T$ is conservative and $R$ is probability preserving,  $T\times R$ is conservative.
Thus we may apply the first claim of the theorem to conclude that $T\times R$ is ergodic.
Hence $T$ is weakly mixing.
\qed
\enddemo

\comment
Let $S$ be a nonsingular one-to-one transformation of $X$ such that
$\{S^nx\mid n\in\Bbb Z\}=\{T_gx\mid g\in G\}$ for a.e. $x$.
It exists according to the main result of \cite{CoFeWe}.
Then by Lemma~2.3, $S$ (considered as a $\Bbb Z$-action) is strictly $\Cal R$-outer, $\Delta_\Cal R$-bounded and squashes $\Cal R$. 
We note that $S$ is conservative because $T$ is conservative.
In a similar way, $S$ is ergodic if and only if $T$ is ergodic.
Thus it remains to prove that $S$ is ergodic.
For that we repeat the proof of Theorem~1.1 almost literally.
The only difference is that we approximate $B$ with a sequence of open subsets $(A_n)_{n=1}^\infty$ for which there is a sequence of positive reals  $(r_n)_{n=1}^\infty$
such that each $\bigcup_{x\in A_n}\{y\in X\mid d(y,x)\le r_n\}=A_n$.
Since $S$ squashes $\Cal R$, it follows that for each $(x,y)\in\Cal R$, we have that  $S^kx\in A_n$ if and only if $S^ky\in A_n$ eventually in $k$ for each $n$.
\endcomment

Given a countable set $A$, we let $X=A^G$ and 
denote by  $\Cal T$  the tail equivalence relation on $X$.
Given a sequence $(\mu_n)_{n=1}^\infty$ of non-degenerated probability measures on $A$, 
we let $\mu:=\bigotimes_{n=1}^\infty\mu_n$. 
We will always assume that  $\mu$ is non-atomic or, equivalently,
$$
\sum_{g\in G}\log(\max_{a\in A}\mu_g(a))=-\infty.
$$
Then
$\Cal T$ is $\mu$-nonsingular and for all $\Cal T$-equivalent points $x=(x_h)_{h\in G}$ and
$y=(y_h)_{y\in H}$ from $X$,
$$
\Delta_{\Cal T,\mu}(x,y)=\prod_{h\in G}\frac{\mu_h(y_h)}{\mu_h(x_h)}.\tag2-5
$$
We note that the product is finite indeed.
It follows from the Kolmogorov 0-1
 law that $\Cal T$ is $\mu$-ergodic\footnote{Hint: enumerate the elements of $G$ with natural numbers.}.
Let $T=(T_g)_{g\in G}$ denote the left shiftwise $G$-action on  $X$.
It follows from \cite{Ka} that $T$ is nonsingular if and only if for each $g\in G$,
$$
\sum_{h\in G}\sum_{a\in A}\Big(\sqrt{\mu_h(a)}-\sqrt{\mu_{g^{-1}h}(a)}\,\Big)^2<\infty.
$$

\comment
Moreover, if $T$ is nonsingular then for each $g\in G$,
$$
\frac{d\mu\circ T_g}{d\mu}(x)=\prod_{h\in G}\frac{\mu_{g^{-1}h}(x_h)}{\mu_{h}(x_h)}.\tag2-6
$$
%The dynamical system $(X,\mu,T)$ is called {\it a nonsingular Bernoulli $G$-action.}

\endcomment

\proclaim{Corollary 2.4} Let $(X,\mu,T)$ be a conservative nonsingular Bernoulli $G$-action with
$A$  finite. 
If 
$
\delta:=\inf_{g\in G}\min_{a\in A}\mu_g(a)>0
$
then $T$ is weakly mixing.
 More generally,
  if the ``diagonal'' $G$-action
 $T\times\cdots\times T(p\text{\ times})$ is conservative for some $p>1$ then  it is weakly mixing.
\endproclaim

\demo{Proof} We note that $X$ is a Cantor (compact, ultrametric) space.
%We define a metric $d$ on $X$ by setting
%$d(x,y):=\sum_{l\in\Bbb Z}2^{-|l|}\delta(x_l,y_l)$, where $\delta(.,.)$ is the Kronecker delta.
%Then $d$ is an ultrametric compatible with the topology on $X$.
%We note that$\Delta_\Cal S(x,y)=\prod_{h\in G}\frac{\mu_h(y_h)}{\mu_h(x_h)}$for all $(x,y)\in\Cal S$.
Of course, $T$ normalizes and  squashes $\Cal T$ (mod $\mu$).
Moreover, $T$ is strictly $\Cal R$-outer (mod $\mu$).
%It remains to verify that $T$ is $\Cal R$-bounded and apply Theorem~2.2.
Take $(x,y)\in\Cal T$.
Let $F_{x,y}:=\{g\in G\mid x_g\ne y_g\}$.
It follows from \thetag{2-5} that 
$$
\Delta_{\Cal T,\mu}(T_gx,T_gy)
%\frac{\frac{d\mu\circ T_g}{d\mu}(x)}{\frac{d\mu\circ T_g}{d\mu}(y)}
=\prod_{h\in G}
\frac{\mu_{gh}(y_{h})}{\mu_{gh}(x_{h})}%\frac{\mu_h(y_h)}{\mu_h(x_h)}
=
\prod_{h\in F_{x,y}}
\frac{\mu_{gh}(y_{h})}{\mu_{gh}(x_{h})}
%\prod_{h\in F_{x,y}}
%\frac{\mu_{g^{-1}h}(x_{h})}{\mu_{g^{-1}h}(y_{h})}\frac{\mu_h(y_h)}{\mu_h(x_h)}
.
$$
Hence
$
\delta^{\#F_{x,y}}
%<\frac{d\mu\circ T_g}{d\mu}(x)/\frac{d\mu\circ T_g}{d\mu}(y)
\le\Delta_{\Cal T,\mu}(T_gx,T_gy)
\le\delta^{-\#F_{x,y}}.
$
Thus, 
 $T$ is $\Cal T$-bounded by Remark~2.2.
It now follows from Theorem~2.3 that $T$ is weakly mixing.
Thus, the first claim is proved.
The second claim follows from the first one.
\qed
\enddemo

\head 3.  On weak mixing of conservative two-sided shifts\endhead

In this section we consider in more detail the case  where $G=\Bbb Z$.
First we note that in this case the condition of the $\Cal R$-boundedness in Theorem~2.3 can be slightly relaxed.

\proclaim{Theorem 3.1} Let $X$ be the same as in Theorem~2.2 and let $\mu$ be a probability Borel measure on $X$.
 Let $T$ be a conservative nonsingular  transformation of $(X,\mu)$ that
 normalizes  a countable  Borel equivalence relation
 $\Cal R$ on $X$ and let the restriction of $T$ to \rom{Inv}$(\Cal R)$ be $\mu$-ergodic.
 If $T$ is strictly $\Cal R$-outer, squashes $\Cal R$  and 
 there is a Borel map $\alpha:\Cal R\to [1,+\infty)$ and a $\mu$-conull Borel subset $X_0\subset X$ such that
$$
\alpha(x,y)^{-1}<\frac{d\mu\circ T^{-k}}{d\mu}(x)/\frac{d\mu\circ T^{-k}}{d\mu}(y)<\alpha(x,y)
%\tag{3-1}
$$
 for all $(x,y)\in\Cal S\cap (X_0\times X_0)$ and each $k\ge 0$
 then $T$ is ergodic.
If, moreover, $\Cal R$ is $\mu$-ergodic then $T$ is weakly mixing.
\endproclaim

\demo{Proof}  The proof is almost a verbal repetition of the proof of Theorem~2.3.
Only %two slight  
a slight modification  is needed.
%First, since $\Cal T$ is hyperfinite itself (as the orbit equivalence relation of an amenable group action), we set $\Cal S:=\Cal T$.
%Second, 
Instead of \thetag{2-2}, we should use the following limit
$$
E(1_{A_n}\mid\goth F)(x)=\lim_{K\to+\infty}\frac{\sum_{k=0}^K1_{A_n}(T^{-k}x)\frac{d\mu\circ T^{-k}}{d\mu}(x)}{\sum_{k=0}^K\frac{d\mu\circ T^{-k}}{d\mu}(x)}
$$
at a.e. $x\in X$ for each $n>0$.
The existence of this limit follows from the  
 Hurewicz ergodic theorem.
 \qed
\enddemo

In a similar way one can ``relax'' the corresponding $\Cal R$-boundedness condition in Remark~2.2.

We now discuss some applications of Theorem~3.1 to  nonsingular Bernoulli shifts.
 The first claim of the following corollary strengthens \cite{Ze, Theorem~3} (where it was proved  that $T$ is ergodic under a stronger condition that $\inf_{n\in\Bbb Z}\min_{a\in A}\mu_n(a)>0$). 

\proclaim{Corollary 3.2} Let $(X,\mu,T)$ be a conservative nonsingular Bernoulli $\Bbb Z$-shift with $A$ finite.
If
$
\delta:=\inf_{n>0}\min_{a\in A}\mu_n(a)>0
$  
then
$T$ is weakly mixing.
More generally, if the product $T\times\cdots\times T (p\text{ times})$ is conservative for some $p>1$ then
 $T\times\cdots\times T (p\text{ times})$ is weakly mixing.
\endproclaim

\demo{Proof} Take $(x,y)\in\Cal T$.
Let $J_{x,y}:=\{k\in\Bbb Z\mid x_k\ne y_k\}$.
It follows from \thetag{2-5} that
$$
%\frac{d\mu\circ T^{-n}}{d\mu}(x)/\frac{d\mu\circ T^{-n}}{d\mu}(y)
\Delta_{\Cal T,\mu}(T^nx,T^ny)=%\prod_{k\in J_{x,y}}\frac{\mu_k(y_k)}{\mu_k(x_k)}
\prod_{k\in J_{x,y}}\frac{\mu_{k+n}(y_k)}{\mu_{k+n}(x_k)}.
$$
If $n$ is sufficiently large then $k+n>0$ for each $k\in J_{x,y}$.
Hence we have eventually (in $n$) that
$
 \delta^{\# J_{x,y}}\le
 \Delta_{\Cal T,\mu}(T^nx,T^ny)
 % \prod_{k\in J_{x,y}}\frac{\mu_{k-n}(x_k)}{\mu_{k-n}(y_k)}
  \le \delta^{-\# J_{x,y}}
$
%and \thetag{3-1} follows.
It remains to apply Theorem~3.1 (and a corresponding analogue of Remark~2.2).
\qed
\enddemo

\comment
We recall that $\mu$ (as well as  the system $(X,\mu,T)$) is called {\it semi-stationary} if there is $n\in\Bbb Z$ such that either $\mu_n=\mu_{n-1}=\mu_{n-2}=\cdots$ or  $\mu_n=\mu_{n+1}=\mu_{n+2}=\cdots$.
A lot of examples of nonsingular Bernoulli shifts not admitting an equivalent invariant $\sigma$-finite measure were constructed in \cite{Ha}, \cite{Ko}, \cite{DaLe}, \cite{VaWa}.
All of them are semistationary. 
\endcomment
We now give some  concrete examples of weakly mixing nonsingular Bernoulli shifts without an equivalent semistationary measure. 
%For that we will apply the conservativeness criterion for Bernoulli shifts from \cite{VaWa}.

\example{Example 3.3}
Fix $\lambda\in(0,1/2)$ and a real $\xi$ such that 
$$
\xi>\left(\frac\lambda 2\right)^{-2}+\left(\frac\lambda 2\right)^{-1}\left(1-\frac\lambda 2\right)^{-2}.
$$
Choose a  sequence of reals $(\alpha_n)_{n=1}^\infty$ and a sequence of positive integers $(A_n)_{n=0}^\infty$ such that the following conditions hold:
\roster
\item"(i)"
$\alpha_j\in (0,1/2)$ and 
$\sum_{j=1}^\infty\alpha_j^2<+\infty$,
\item"(ii)"  $A_0=0$, $A_l>8A_{l-1}$ for each $l>0$ and 
 %, $A_n>4$ for each $n$ and $\sum_{n=1}^\infty\frac1{A_n}<\infty$.
\item"(iii)" %$1=\alpha_1^2A_1=\alpha_2^2A_2=\cdots$.
$\sum_{j=1}^n\alpha^2_jA_j=\frac{\log n}{4\xi}$ for each $n>0$.
\endroster
We now define for each $n\in\Bbb Z$, a measure $\mu_n$ on $A=\{0,1\}$ by setting
$$
\mu_n(0):=
\cases
\lambda,&\text{if $2A_l\le |n|<A_{l+1}$ for some $l\ge 0$},\\ 
\lambda+\alpha_l, &\text{if $A_l\le |n|<2A_l$ for some $l> 0$.}
\endcases
$$
It is straightforward to verify that \roster
\item"(I)"
$\lambda/2<\mu_n(0)<1-\lambda/2$ for each $n\in\Bbb Z$
and
\item"(II)"
$\sum_{n\in\Bbb Z}|\mu_{n+1}(0)-\mu_n(0)|^2=4\sum_{j=1}^\infty\alpha_j^2<\infty$.
\endroster
It follows from (I) that $\mu$ is non-atomic.
It follows from \thetag{2-5} and (II) that $T$ is $\mu$-nonsingular.
Of course,  $\mu$ is not semi-stationary.
It is straightforward to deduce from (I), (iii) and the Kakutani theorem \cite{Ka}  that $\mu$ is not equivalent to any semi-stationary measure on $X$.
To show that $T$ is conservative we will apply \cite{VaWa, Proposition~4.1}, which
 states that if $\sum_{n\in\Bbb Z}e^{-\xi\sum_{j\in\Bbb Z}|\mu_{j+n}(0)-\mu_j(0)|^2}=+\infty$
 then $T$ is conservative.
Indeed, it follows from (ii) and (iii) that for each $l>0$,
$$
\sum_{j\in\Bbb Z}|\mu_{j-4A_l}(0)-\mu_j(0)|^2=4\sum_{j=1}^l\alpha_j^2A_j+4\sum_{j>l}\alpha^2_jA_l=\frac{\log l}{\xi}+\overline o(1).
$$
Then
$$
\align
\sum_{n\in\Bbb Z}e^{-\xi\sum_{j\in\Bbb Z}|\mu_{j+n}(0)-\mu_j(0)|^2}&\ge
\sum_{l=1}^{+\infty}e^{-\xi\sum_{j\in\Bbb Z}|\mu_{j-4A_{l}}(0)-\mu_j(0)|^2}\\
&=\sum_{l=1}^{+\infty}e^{-\log l+\overline o(1) }
\endalign
$$
Hence
$
\sum_{n\in\Bbb Z}e^{-\xi\sum_{j\in\Bbb Z}|\mu_{j+n}(0)-\mu_j(0)|^2}=+\infty
$
and 
thus $T$ is $\mu$-conservative.
Corollary~3.2 yields that $T$ is weakly mixing.
\endexample

We now investigate ergodicity of nonsingular Bernoulli shifts under
assumption that $\delta=0$, where $\delta$ is defined in  Corollary~3.2.

Let $\Cal S$ be {\it the symmetric} equivalence relation on $X$.
By definition, it  is the orbit equivalence relation of the natural action of the  group $\Sigma_0({\Bbb Z})$ of
all finite permutations of $\Bbb Z$ on $X$.
In other words, two points $x=(x_n)_{n\in\Bbb Z}$ and
$y=(y_n)_{n\in\Bbb Z}$ of $X$ are $\Cal S$-equivalent if there is $\sigma\in \Sigma_0({\Bbb Z})$
such that $x_n=y_{\sigma(n)}$ for all $n\in\Bbb Z$.
It is easy to verify that $\Cal S$ is  a subrelation of $\Cal T$ and $T$ normalizes $\Cal S$.
%We consider first the simplest case where $A$ consists of two points only.
We will use the following result from \cite{AlPi}.

\proclaim{Lemma 3.4} 
\roster
\item"\rom(i)" Let $\# A=2$.
Then $\Cal S$ is $\mu$-ergodic if and only if 
$
\sum_{n\in\Bbb Z}\min_{a\in A}\mu_n(a)=+\infty,
$
i.e. $\mu$ is non-atomic.
\item"\rom(ii)" Let $A$ be finite.
Then $\Cal S$ is $\mu$-ergodic if and only if for each nonempty proper subset $B\subset A$, 
$
\sum_{n\in\Bbb Z}\min(\mu_n(B),\mu_n(A\setminus B)=+\infty.$
\endroster
\endproclaim

The following statement strengthens   Corollary~3.2 in the case where $\# A=2$.

\proclaim{Theorem 3.5} Let $(X,\mu, T)$ be a nonsingular Bernoulli shift  with $\# A=2$. 
If  there is $D>1$ such that
$$
D^{-1}<\frac{\mu_k(a)}{\mu_{k+1}(a)}<D\quad\text{for all $a\in A$ and each }k<0
\tag3-1
$$
and $T$ is conservative 
  then
$T$ is weakly mixing.
More generally, if $T\times\cdots\times T (p\text{ times})$ is conservative for some $p>1$ then
 $T\times\cdots\times T (p\text{ times})$ is weakly mixing.
\endproclaim
\demo{Proof}
%We first note that $\sum_{n\in\Bbb Z}\min_{a\in A}\mu_n(a)=+\infty$.
%Indeed, otherwise we would have that  $\infty>\sum_{n\in\Bbb Z}\log(1-\min\mu_n(a))=
%\sum_{n\in\Bbb Z}\log(\max\mu_n(a))$.
%This implies that $\mu$ is atomic, a contradiction.
%Therefore, 
By Lemma~3.4(i), $\Cal S$ is ergodic.
Take $(x,y)\in\Cal S$.
Let $J_{x,y}:=\{k\in\Bbb Z\mid x_k\ne y_k\}$ and let $\sigma$ be a permutation of $J_{x,y}$
such that $y_k=x_{\sigma(k)}$ for all $k\in J_{x,y}$.
It follows from \thetag{2-5} that
$$
%\frac{d\mu\circ T^{-n}}{d\mu}(x)/\frac{d\mu\circ T^{-n}}{d\mu}(y)
\Delta_{\Cal S,\mu}(T^nx,T^ny)=%\prod_{k\in J_{x,y}}\frac{\mu_k(y_k)}{\mu_k(x_k)}
\prod_{k\in J_{x,y}}\frac{\mu_{\sigma^{-1}(k)+n}(x_k)}{\mu_{k+n}(x_k)}.
$$
If $n$ is large enough, we deduce from \thetag{3-1} that
$$
D^{-|\sigma^{-1}(k)-k|}\le\frac{\mu_{\sigma^{-1}(k)+n}(x_k)}{\mu_{k+n}(x_k)}\le D^{|\sigma^{-1}(k)-k|}
$$
for each $k\in J_{x,y}$.
Since $|\sigma^{-1}(k)-k|<\# J_{x,y}$, we obtain that
$$
D^{-(\# J_{x,y})^2}\le \frac{\mu_{\sigma^{-1}(k)+n}(x_k)}{\mu_{k+n}(x_k)}\le D^{(\# J_{x,y})^2}
$$
eventually in $n$.
Hence  $T$ is  weakly mixing by Theorem~3.1 and the analogue of Remark~2.2.
\qed
\enddemo

As in  \cite{Ko4}, we note that Theorem~3.1 can be applied also to nonsingular Markov shifts.

%\head 4. Nonsingular Markov shifts
%\endhead

\remark{Remark \rom{3.6}}
Let $ S$ be a finite set and let $M=(M(a,b))_{a,b\in S}$ be a $\{0,1\}$-valued $(S\times S)$-matrix.
Given two integers $i\le j$ and a finite sequence $a=(a_l)_{i\le l\le j}$ of elements from $S$ such that $M(a_l,a_{l+1})=1$ for each $l=i,\dots,j-1$, we denote by
$[a]_i^j$ the cylinder $[a]_i^j:=\{x\in X_M\mid x_l=a_l\text{ for each }l=i,\dots,j\}$.
The set  of infinite paths $X_M:=\{x=(x_i)_{i\in\Bbb Z}\in S^\Bbb Z\mid M(x_i,x_{i+1})>0\}$
 is a closed subset of $S^\Bbb Z$ which is invariant under the two-sided shift $T$.
 %It is invariant under the two-sided shift $T$ an $S^\Bbb Z$.
%The topological dynamical system $(X_M,T)$ is called a {\it shift of finite type (SFT)}.
%Given two integers $i\le j$ and a finite sequence $a=(a_l)_{i\le l\le j}$ of elements from $S$ such that $M(a_l,a_{l+1})>0$ for each $l=i,\dots,j-1$, we denote by
%$[a]_i^j$ the cylinder $[a]_i^j:=\{x\in X_M\mid x_l=a_l\text{ for each }l=i,\dots,j\}$.
%The set of all cylinders is a base of the topology on $X_M$.
%This implies that  the corresponding  SFT is topologically mixing.
Suppose  that there is  a sequence $(\pi_n)_{n\in\Bbb Z}$ of probability measures on $S$
and a sequence $(P_n)_{n\in\Bbb Z}$ of row-stochastic $(S\times S)$-matrices such that $\pi_nP_n=\pi_{n+1}$ 
and  $P_n(a,b)>0$ if and only if $M(a,b)=1$ for each $n\in\Bbb Z$.
Then there is a unique probability measure $\mu$ on $X_M$ such that for every cylinder
$[a]_i^j$ in $X_M$,
$$
\mu([a]_i^j)=\pi_i(a_i)P_i(a_{i},a_{i+1})\cdots P_{j-1}(a_{j-1},a_j).
$$
%Of course, $\mu$ is fully supported on $X_M$.
It is called   {\it a Markov measure} on $X_M$ generated by $(\pi_n,P_n)_{n\in\Bbb Z}$.
Denote by $\Cal T$ the restriction of the tail equivalence relation to $X_M$.
It is easy to verify that $\Cal T$ is $\mu$-nonsingular and
for each $(x,y)\in\Cal T$,
$$
\Delta_{\Cal T,\mu}(x,y):=\prod_{i\in\Bbb Z}\frac{P_i(y_i,y_{i+1})}{P_i(x_i,x_{i+1})}.
$$
We note that the infinite product in this formula is indeed finite.
Suppose that 
\roster
\item"---"
$M$ is {\it primitive}, i.e.   there is $n>0$ with
$M^n(a,b)>0$ for all $a,b\in S$ and
\item"---" that $\delta:=\inf \{P_n(a,b)\mid n\in\Bbb Z, M(a,b)=1\}>0$.
\endroster
It was shown in \cite{Ko4} that under these two conditions, $\Cal T$ is $\mu$-ergodic. 
It follows from this fact, Theorem~3.1 and Remark~2.2 that if  the shift $T$ on $(X_M,\mu)$ is conservative and  nonsingular then $T$ weakly mixing.
%Let $\Cal S$ stand for {\it the symmetric} equivalence subrelation of $\Cal S$.
%By definition, two points $x=(x_i)_{i\in \Bbb Z}$ and $y=(y_i)_{i\in \Bbb Z}$ from $X_M$ are $\Cal S$-equivalent if there is a  finite permutation $\sigma$ of $\Bbb Z$ such that $x_i=y_{\sigma(i)}$ for each $i\in\Bbb Z$.
Indeed, it suffices to note that 
$\Delta_{\Cal T,\mu}(T^nx,T^ny):=\prod_i\frac{P_{i+n}(y_i,y_{i+1})}{P_{i+n}(x_i,x_{i+1})}$
and argue as in Corollary~3.2.
\endremark

\comment
\proclaim{Theorem 4.1} $\Cal S$ is $\mu$-ergodic.
\endproclaim

\demo{Proof} Fix  $s_1,s_2\in S$ such that $M(s_1,s_2)>0$.
Given $n\ge 0$ and two cylinders $[a]_{-n}^n$ and $[b]_{-n}^n$.
%Since $M^l>0$ for some $l>0$, 
there are cylinders $[a']_{-n-l}^{n+l}$ and
$[b']_{-n-l}^{n+l}$ such that $a'_i=a_i$ and $b'_i=b_i$ if $|i|\le n$, 
$$
a'_{n+l}=b'_{n+l}=s_1\ \text{ and } \ 
a'_{-n-l}=b'_{-n-l}=s_2.
$$

\enddemo

Under some conditions on $(\pi_n,P_n)_{n\in\Bbb Z}$, the corresponding Markov measure is nonsingular.
We then call the dynamical system $(X_M,\mu,T)$ {\it the nonsingular Markov shift}.

\proclaim{Theorem 4.1}  If a nonsingular Markov shift is conservative then it is ergodic.
If $M$ is primitive (i.e. the corresponding SFT is topologically mixing) then $T$ is weakly mixing.
\endproclaim
\demo{Proof}
\enddemo

\endcomment

\head Appendix. Weak pointwise ratio ergodic theorem for \\ nonsingular actions of amenable  groups 
\endhead

%We first prove a weak version of  the pointwise ergodic theorem for the nonsingular actions
%of amenable groups. 
We recall that given a finite subset $K\subset G$ and $\epsilon>0$, a finite subset $F\subset G$ is called {\it $[K,\epsilon]$-invariant} if $\#\{g\in F\mid Kg\subset F\}>(1-\epsilon)\#F$.

\proclaim{Theorem A.1} Let $G$ be amenable.
Fix a sequence $(K_n)_{n=1}^\infty$ of finite subsets in $G$ and a sequence of positive reals
$(\epsilon_n)_{n=1}^\infty$ converging to $0$.
Then for each  nonsingular action $T=(T_g)_{g\in G}$ on a standard probability space $(X,\mu)$ and a countable subset $\Cal L$ of $L^1(X,\mu)$, there is an increasing sequence $F_1\subset F_2\subset\cdots$ of finite subsets in $G$ such that $\bigcup_{n=1}^\infty F_n=G$ and for every $f\in\Cal L$, there exists
$$
\lim_{n\to\infty}\frac{\sum_{g\in F_n}f( T_gx)\frac{d\mu\circ T_g}{d\mu}(x)}{\sum_{g\in F_n}\frac{d\mu\circ T_g}{d\mu}(x)}=E(f\mid\text{\rom{Inv}}(\Cal T))(x),\tag A-1
$$
at a.e. $x$, where $\Cal T$ denotes the $T$-orbit equivalence relation on $X$.
Moreover, $F_n$ is $[K_m,\epsilon_m]$-invariant for each $m\le n$.
In particularly, if $\bigcup_{n=1}^\infty K_n=G$ then $(F_n)_{n=1}^\infty$ is a F{\o}lner sequence in $G$.
\endproclaim

\demo{Proof}
Let $R=(R_g)_{g\in G}$ be an ergodic measure preserving Bernoulli (free) action of $G$ on a standard  probability space $(Y,\nu)$.
Denote by $\Cal R$ the $R$-orbit equivalence relation on $Y$.
We also denote by $\widehat{\Cal R}$ the orbit equivalence relation of the product $G$-action
$T\times R:=(T_g\times R_g)_{g\in G}$.
Since $G$ is amenable, $\Cal R$ and $\widehat{\Cal R}$ are both hyperfinite.
Let $(\Cal R_n)_{n=1}^\infty$ be a filtration of $\Cal R$. 
For $y\in Y$ and $n>0$, let  $G_{y,n}$ stand for  the finite subset of $G$ such that
$\Cal R_n(y)=\{R_gy\mid g\in G_{y,n}\}$.
For each $n>0$, we now define a finite  equivalence relation $\widehat {\Cal R}_n$ on $X\times Y$ by setting
$$
\multline
(x,y)\sim_{\widehat{\Cal R}_n}(x',y')
\quad\text{iff  $(y,y')\in\Cal R_n$
 and $x'=T_gx$} \\
\text{ for the
 only element $g\in  G$ such that
  $y'=T_gy$. }
  \endmultline
$$
Then $(\widehat{\Cal R}_n)_{n=1}^\infty$ is  a filtration of  $\widehat{\Cal R}$.
It follows from the martingale convergence theorem there there is a $(\mu\times\nu)$-conull subset $Z$ of $X\times Y$ such that
$$
\lim_{n\to\infty}E(f\otimes 1\mid\text{Inv}(\widehat{\Cal R}_n))(x,y)= E(f\otimes 1\mid\text{Inv}(\widehat{\Cal R}))(x,y)
\tag A-2
$$
for every $f\in\Cal L$ at each $(x,y)\in Z$.
We now note that
$$
\aligned 
E(f\otimes 1\mid\text{Inv}(\widehat{\Cal R}_n))(x,y)&=\frac{\sum_{g\in G_{y,n}} f(T_gx)\frac{d\mu\circ T_g}{d\mu}(x)}{\sum_{g\in G_{y,n}}\frac{d\mu\circ T_g}{d\mu}(x)}\quad\text{and}\\
E(f\otimes 1\mid\text{Inv}(\widehat{\Cal R}))(x,y)&=E(f\mid\text{Inv}(\Cal T))(x).
\endaligned
\tag A-3
$$
By the Fubini theorem, there is $y_0\in Y$ and a $\mu$-conull subset $X_0$ in $X$
such that  $X_0\times\{y_0\}\subset Z$. 
We let $F_n:=G_{y_0,n}$.
Since $\Cal R_1(y_0)\subset\Cal R_2(y_0)\subset\cdots$ and $\bigcup_{n=1}^\infty\Cal R_n(y_0)=\Cal R(y_0)$, it follows that $F_1\subset F_2\subset\cdots$ and $\bigcup_{n=1}^\infty F_n=G$.
It remains to substitute ~\thetag{A-3}  into \thetag{A-2}  to obtain \thetag{A-1}.

By \cite{Da, Lemma~2.2}, for each $m>0$, 
$$
\nu(\{y\in Y\mid G_{y,n}\text{ is $[K_m,\epsilon_m]$-invariant}\})>1-\epsilon_m
$$
eventually in $n$.
Hence passing to a subsequence in $(\Cal R_n)_{n=1}^\infty$ and applying the Borel-Cantelli lemma we can choose $y_0$ in such a way that  $F_n$ is $[K_m,\epsilon_m]$-invariant for all $n\ge m>0$.
\qed 
\enddemo

\comment

Now we consider a more general case where $G$ is assumed to be an arbitrary infinite countable group without Kazhdan property $(T)$.

\proclaim{Theorem A2} Let $G$ do not have property $(T)$.
Then for each  nonsingular action $T=(T_g)_{g\in G}$ on a standard probability space $(X,\mu)$ and a countable subset $\Cal L$ of $L^1(X,\mu)$,
there is an increasing sequence $F_1\subset F_2\subset\cdots$ of finite subsets in $G$ such that \thetag{A-1} holds for each $f\in \Cal L$.
\endproclaim

\demo{Proof} Since  $G$ does not have property $(T)$, it follows from \cite{CoWe} that there exists an ergodic probability preserving $G$-action $R=(R_g)_{g\in G}$ which is not strongly ergodic.
Let $(Y,\nu)$ be the space of $R$.
Denote by $\Cal R$ the $R$-orbit equivalence relation on $Y$.
Since $R$ is not strongly ergodic, there is a standard probability space $(Z,\kappa)$, 

\enddemo

\endcomment

\enddocument